\documentclass{elsarticle}
\makeatletter
\def\ps@pprintTitle{%
 \let\@oddhead\@empty
 \let\@evenhead\@empty
 \def\@oddfoot{}%
 \let\@evenfoot\@oddfoot}
\makeatother
\usepackage[margin=2.5cm]{geometry}
\usepackage{lineno,hyperref}
\usepackage[fleqn]{amsmath}
\usepackage{caption}
\usepackage{subcaption}
\usepackage[labelfont=bf]{caption}
\modulolinenumbers[5]
\usepackage{tabularx,booktabs}
\usepackage{graphicx}
\usepackage{subcaption}
\usepackage{amsmath,amsthm,amssymb}
\usepackage{epstopdf}
\usepackage{mathrsfs}
\usepackage{setspace}

\biboptions{authoryear}

\begin{document}
\begin{frontmatter}
\title{A probability-based multi-path alternative fueling station location model}
\author[a]{\corref{mycorrespondingauthor} Shengyin Li}
\cortext[mycorrespondingauthor]{Corresponding author}
\ead{shengyinli@gmail.com}
\author[a]{Yongxi Huang}
\address[a]{Glenn Department of Civil Engineering, 109 Lowry Hall, Clemson University, Clemson, SC 29634, USA}

\begin{abstract}
We develop a probability-based multi-path location model for an optimal deployment of alternative fueling stations (AFSs) on a transportation network. Distinct from prior research efforts in AFS problems, in which all demands are deemed as given and fixed, this study takes into account that not every node on the network will be equally probable as a demand node. We explicitly integrates the probability into the model based on the multi-path refueling location model to determine the optimal station locations with the goal to maximize the expected total coverage of all demand nodes on a system level. The resulting mixed integer linear program (MILP) is NP-hard. A heuristic based on genetic algorithm is developed to overcome the computational challenges. We have conducted extensive numerical experiments based on the benchmark Sioux Falls network to justify the applicability of the proposed model and heuristic solution methods. \\
\end{abstract}
\begin{keyword}
\texttt{location, probability, multi-path, AFS, genetic algorithm}
\end{keyword}
\end{frontmatter}

\section{Introduction}
\noindent The lack of sufficient public fueling stations for Alternative Fuel Vehicles (AFVs) has greatly hindered the adoption of AFVs as mass-market modes of personal transportation. Take electric vehicles (EVs) as an example. Most major auto manufacturers build vehicles with some type of battery propulsion system, that are either hybrid-electric powered or totally battery powered. It is expected that nearly 6 million EVs will be in use worldwide by 2020, of which about 1 million will be in the U.S. \citep{trigg2013etal}. Compared to conventionally fueled vehicles with ranges well in excess of 300 miles, the average vehicle range on EVs is low (e.g., an average of 82 miles on Nissan Leaf ). The success of mass adoptions of EVs thereby heavily relies on the availability of public charging, which currently is woefully inadequate. An effective system of public fueling stations on transportation networks is crucial to facilitate successful mass adoption of AFVs.

On the other hand, the resource of alternative fueling stations (AFSs) is limited and thus strategic locations of AFS is of great research interest recently \citep{hodgson1990flow, kuby2005flow, wang2010locating, huang2015optimal}. Regardless the objective to either maximize the coverage of trips completed by the AFVs or minimize the total cost of fulfilling all trip demands by AFVs, the underlined assumption is that all trip demands are known a priori, by which the origin-destination (O-D) pairs of the trips are fixed and deterministic. However, the trips or O-D pairs, in most cases, are estimated statistically and the resulting probabilities of a node becoming a demand node of interest should be implicative to the decisions of locations of AFSs, especially when resources are limited. In this study, we develop a probability-based optimization model that  rolls out AFSs on a network to maximize the expected coverage of nodes, from which round-trips are able to be completed by AFVs.

Flow-based station location problems originate from node-based facility location problems \citep{owen1998strategic} and treat the demand as flows of goods or services. 
Many models have been developed for use in locating charging stations, the most notable exemplars of which are the Flow Intercepting Facility Location Problem (FIFLP) \citep{hodgson1990flow, berman1992optimal, boccia2009flow}, the maximal-flow-coverage based Flow-Refueling Location Problem (FRLP) \citep{kuby2005flow, kim2012deviation, kim2013network, upchurch2009model}, and the minimum-cost flow based set-covering \citep{wang2007optimal, wang2009locating, wang2008locating, wang2013locating, wang2010locating}, with applications in the Alternative Fueling Station (AFS) problem \citep{kuby2009optimization, lim2010heuristic, zeng2010generalized, capar2012efficient}. Both the maximum-coverage and minimum-cost models were recently generalized \citep{mirhassani2012flexible, wen2013locating}. 
More recently, a series of models have been developed  to enable the integrated decisions on strategic charging locations and travelers' route choices while considering the spatial distributions of demands \citep{huang2015optimal, li2014heuristic} and extended the model to a multi-stage optimization model \citep{Lietal2016}. 

However, none of those models takes into account the effects of uncertainty in the process of planning out AFSs on a network, which could lead to potentially biased strategies. Taking the EV market as an example, when decision makers assess if a city would become a next EV market, the resulting probability could be a multivariate statistical analysis (e.g., logistic regression analysis) and be result of many factors. Large investments on cities with low probability to become the next EV market should be avoided. Thus, probabilistic facility location models are relevant to this study.

Facility location decisions without considering those uncertainties of network characteristics may result in a waste of resources. Among the factors that influence the public facility location decisions, the most apparent one is the stochastic nature of demand. When uncertainty effects are taken into account, a probabilistic approach would provide a more cost-effective system design. In this case, a probabilistic modeling of the system is required in the development of objectives or constraints of the optimization model. In the research community of facility location, there are generally two schools of probabilistic models \citep{marianov20024}. The first is in maximizing expected coverage of each demand node. The second, in either constraining the probability of at least one server being available to be greater than or equal to a specified level $\alpha$, or to count demand nodes covered if this probability is at least $\alpha$. The first category of models were initiated by the maximization of expected covering location problem (MEXCLP) by \citep{daskin1983maximum}, which maximized the expected value of population coverage given limited facilities to be located on the network. Instead of maximizing expected value of coverage, the probabilistic location set covering problem (PLSCP) \citep{revelle1988reliability} was a pilot study in the second category, in which the probability that at least one serving being available to each demand node was constrained to be greater than or equal to a reliability level $\alpha$. These probabilistic models have been applied to different application domains. For example, \cite{belardo1984partial,  psaraftis1986optimal} incorporated a relative likelihood of spill occurrence in a partial covering approach to site response resources such that the probability to cover a spill event was maximized. Readers are referred to two excellent survey papers \citep{daskin2005facility, snyder2006facility} for studies on facility location under uncertainty. More recently, there is a new school of reliable facility location models introduced to relax the assumption that all the facilities are always available.  \cite{snyder2005reliability} introduced the reliable fixed-charge location problem (RFLP), in which facilities are subject to random disruptions. This model has been extended to the recent variants \citep{cui2010reliable, li2010continuum, shen2011reliable}. See \citep{snyder2015or} for a review of the literature on supply chain management with disruptions.

 

In this study, we develop a probability-based alternative fueling station location model. The study is motivated by the fact that the decision on locations of AFSs on the network are driven by the prediction of trip demands, which the results of statistical analysis of demographic and economic data. To our knowledge, this is the first model of this kind by adding probabilistic concern to the flow-based facility location problem. In particular, we will extend our the multi-path refueling location model (MPRLM) \citep{huang2015optimal} to the probabilistic MPRLM or P-MPRLM, by integrating the probability of each city to become an EV market (or called EV adopter) in an optimization model to strategically locate AFS on the network. The goal is to maximize the total expected coverage of all demand nodes on the network with a given budget (i.e., the maximum number of stations available to be located) while determining the best locations of the AFSs on the network. In this model, a node coverage is defined as the number of round trips  originated from the node to be completed. A completion of a round trip depends on a number of factors: topological structure of the network, willingness to take deviation paths, vehicle range, and more importantly the availability and distribution of AFSs. All these factors are integrated in the proposed modeling framework. The model will be a mixed integer linear program and NP-hard. We then develop a heuristic based on genetic algorithm and conduct extensive numerical experiments based on the benchmark Sioux Falls network, through which we explore the interplay between geographic distributions of cities, vehicle range, deviation choice, and probability of becoming EV market. 

The remainder of the paper is organized as follows. In section 2, the formulation of the P-MPRLM is presented and the heuristic based on genetic algorithm is described in section 3. In section 4, setup and results of numerical experiments on Sioux Falls network are reported. We conclude the study and briefly outline the directions of future work in section 5.
\section{Mathematical formulation}
\noindent The new, P-MPRLM is extended from the MPRLM \citep{huang2015optimal} in two major ways. First, the  model is in recognition of the probability of a node (e.g., city) becoming a demand node, which could be a result of multivariate socioeconomic factors, such as income, education level, household size, home ownership, and population density. With limited budgetary resources, it is not economically realistic to satisfy trips between every pair of nodes. Second, round-trips between node pairs are considered in this  model to ensure all travelers can return to their origins. In addition, we also consider and integrate deviation paths into the definition of coverage of a node.  The model aims to maximize the expected coverage of the nodes, from which travelers can freely travel with alternative fueling vehicles (e.g., electric vehicles). In the following subsections, we first define the coverage of a node (or city) on network in section 2.1, followed by the the P-MPRLM formulation in section 2.2. 


\subsection{Definition of node expected coverage}
\begin{figure}[ht]
\centering
\includegraphics[width=0.4\textwidth]{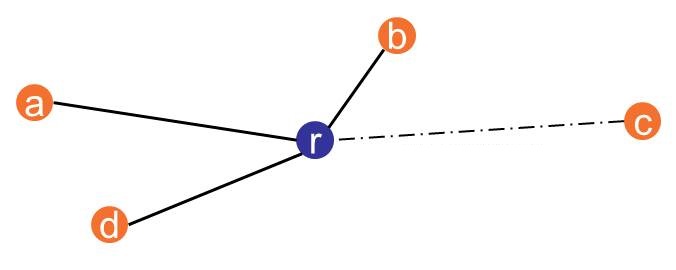}
\caption{An example of a node coverage}
\label{DefinitionOfCoverage}
\end{figure}

\noindent 
In this study, a trip between an O-D pair can be completed via shortest and deviation paths with or without refueling and this O-D pair is considered as ``covered". The \emph{coverage of a node}, say $r$, is then defined as the percentage of covered O-D pairs originated from the node $r$.   
This definition is illustrated by using a small network in Figure \ref{DefinitionOfCoverage}, which consists of a origin node $r$ and four destination nodes $a$, $b$, $c$, and $d$. We use solid lines to indicate that there exists \textit{at least} one path (shortest or deviation paths) that can be used to complete a round-trip between node $r$ and destination nodes, (i.e., O-D pairs $r-a$, $r-b$, and $r-d$), and a dashed line to denote a case, in which trips between the O-D pair cannot be completed, i.e., ($r-c$). As a result, the coverage of node $r$ is $3/4$, denoted as $z_r=3/4$. A node is fully covered only if drivers from the node can complete round-trips to all other nodes on the network with or without refueling. If $p_r$ is the probability of node $r$ becoming a demand node, the expected coverage of node is defined as $p_rz_r$ or $p_r\times3/4$ in this particular case. The details on finding deviation paths and the integrated modeling of deviation paths in location decisions can be referred to the  recent study \citep{huang2015optimal}.

\subsection{P-MPRLM formulation}
\noindent 
The model aims to maximize the expected coverage of nodes on the network through strategically locating limited number of AFSs to facilitate the trips by vehicles with limited driving ranges. It is essentially a maximum-coverage problem, subject to vehicle range, deviation tolerance (i.e., maximum deviation from the planned route), and budget constraints. 
The model is constructed based on a few assumptions. First, a round-trip will be completed by using the same path. This assumption is justifiable as drivers tend to use the paths that they are more familiar with and holds for inter-city and intra- or inter- state trips, where in most cases, driver would complete the trips by using the same inter-state highway systems. Second, we assume that vehicles are homogeneous in terms of the vehicle range and initial state-of-refueling (SOR) at origins. We later run sensitivity analysis on different initial SORs and see how this impact the model results. Third, energy consumed is unified in terms of travel distance. Fourth, this study is focused on the spatial relationship between O-D pairs and roadway networks or as known as spatial economics study, traffic flows are not taken into account, nor are the refueling and travel costs of paths. Extension to explicitly considering such costs may be fulfilled in the future when individual trips can be traced with the aid of large-scale data, such as call detailed records (CRDs). As in this study, we thereby assume that all AFSs are uncapacitated.

Let ($N, A$) be a transportation network, where $N$ and $A$ are the sets of nodes and links, respectively. Let $\tilde{N}$ be the set of candidate station locations, $\tilde{N} \subseteq N$. Let $R\subseteq N$, index $r$, be a set of origin nodes, and $S\subseteq N$, index $s$, be a set of destination nodes. Let $K^{rs}$ be a predefined maximum number of deviation paths for O-D pair $r-s$, which are exogenously generated \citep{huang2015optimal}. We denote by $P^{rs,k}$ a sequence of nodes on the $k^{th}$ path for O-D pair $r-s$, where $k = 1, 2,…, K^{rs}$. Denote a link by $a\in A$ or a pair of ending nodes, i.e., $a=(i,j) \in A$.

The notation used in the model is first presented as follows. We then present the complete mathematical formulation of the P-MPRLM in (1)-(15).  \\

\noindent \emph{Indices:}\\
$i$: index of candidate sites, $i\in \tilde{N} \subseteq N$,\\
$r$: index of origin in the network, $r \in R \subseteq N$,\\
$s$: index of destination in the network, $s \in S \subset N$,\\
$k$: index of the paths for an O-D pair, $k=1,2,      \ldots, K^{rs}$,\\
$a$: index of arc set A, $a=(i,j) \in A$,\\

\noindent \emph{Parameters}:\\
$C_i$: the installing cost of a refueling station, $i \in \tilde{N}$,\\
$\beta$: onboard fuel capacity (unified in travel distance), i.e., vehicle range,\\
$\beta_0$: initial state-of-refueling,\\
$m$: budget,\\
$p_r$: probability of a node to be demand node $r$,\\
$n_r$: predefined number of destinations for demand node $r$,\\
$P^{rs,k}$: a sequence of nodes on the $k^{th}$ path from $r$ to $s$ and then back to $r$ by the same path, where $k=1,2, \ldots, K^{rs}$,\\
$d_{ij}$: distance between node $i$ and $j$,\\
$\delta_i^{rs,k}$: =1 if node $i$ is in the set of node $P^{rs,k}$, 0 otherwise; this is an outcome of the deviation paths that are exogenously generated,\\

\noindent \emph{Variables}:\\
$X_i$: =1 if an AFS is located at node $i$; 0 otherwise,\\
$Y^{rs,k}$: =1 if the $k^{th}$ path between $r$ and $s$ can be completed or used; 0 otherwise,\\
$y^{rs}$: =1 if O-D pair $r-s$ is covered; 0 otherwise,\\
$z_r$: the coverage of node $r$,\\
$B_i^{rs,k}$: remaining onboard fuel at node $i$ on the $k^{th}$ path of O-D pair $r-s$,\\
$l_i^{rs,k}$: amount refueled at node $i$ on the $k^{th}$ path of O-D pair $r-s$.\

\begin{align}
 & \text{Maximize} 
   ~\sum\limits_{r \in R} p_rz_r
\end{align}
\begin{align}
& \text{Subject to} \nonumber \\
& B_i^{rs,k}+l_i^{rs,k} \leq M(1-Y^{rs,k})+\beta, \    \forall r \in R,s \in S;i\in P^{rs,k}, k=1,2, \ldots, K^{rs}\\
& B_i^{rs,k} + l_i^{rs,k} - d_{ij} - B_j^{rs,k} \leq M(1-Y^{rs,k}) ,  \    \forall r \in R,s \in S; i,j \in P^{rs,k}; (i,j)\in A, k=1,2,\ldots, K^{rs}\\
& -(B_i^{rs,k} + l_i^{rs,k} - d_{ij} - B_j^{rs,k}) \leq M(1-Y^{rs,k}), \    \forall r \in R,s \in S; i,j \in P^{rs,k}; (i,j)\in A, k=1,2, \ldots, K^{rs}\\
& \sum\limits_{r \in R} \sum\limits_{s \in S} \sum\limits_{k=1}^{K^{rs}} l_i^{rs,k} \delta_i^{rs,k} \leq MX_i, \    \forall i \in \tilde{N}\\
& \sum\limits_{k=1}^{K^{rs}}Y^{rs,k} \leq My^{rs}, \    \forall r \in R,s \in S\\
& y^{rs} \leq \sum\limits_{k=1}^{K^{rs}}Y^{rs,k},  \ \forall r \in R,s \in S\\
& z_r = \frac{\sum\limits_{s \in S}y^{rs}}{n_r}, \    \forall r \in R\\
& \sum\limits_{i \in \tilde{N}}C_iX_i \leq m \\
& B_i^{rs,k}= \beta_0, \    \forall r \in R,s \in S; i \in R; k=1,2, \ldots, K^{rs}\\
& X_i=\{0,1\}, \    \forall i \in \tilde{N}\\
& Y^{rs,k}=\{0,1\}, \    \forall r \in R,s \in S; k=1,2,\ldots,K^{rs}\\
& y^{rs}=\{0,1\}, \    \forall r \in R,s \in S\\
& z_r \in [0,1], \forall r \in R\\
& B_i^{rs,k} \geq 0, l_i^{rs,k} \geq 0, \    \forall r \in R,s \in S; i \in P^{rs,k} 
\end{align}
\\
\indent The objective of the model is to maximize the expected coverage of demand nodes. Constraint set (2) assures that the onboard fuel of each vehicle does not exceed the fuel capacity (i.e., $B_i^{rs,k}+l_i^{rs,k} \leq \beta$) on the used path ($Y^{rs,k}=1$); otherwise no restriction applies when $Y^{rs,k}=0$. Constraints (3) and (4) work simultaneously to ensure that the fuel consumption conservation $B_i^{rs,k}+l_i^{rs,k}-d_{ij}-B_j^{rs,k}=0$ holds for all links traversed on the $k^{th}$ path which is taken to deploy adequate stations ($Y^{rs,k}=1$); otherwise, if $Y^{rs,k}=0$, then $B_i^{rs,k}+l_i^{rs,k}-d_{ij}-B_j^{rs,k} \leq M$, namely no restraining effects. Constraint set (5) is a logic constraint, stating that refueling only takes place at node $i$ if there is an AFS open there. Constraints (6) and (7) establish the relationship between $Y^{rs,k}$ and $y^{rs}$. In particular, as long as there is a (deviation) path existing between an O-D pair, by which trips between this O-D pair can be completed ($Y^{rs,k}=1$), this O-D pair is deemed as ``covered" ($y^{rs}=1$); otherwise, the O-D pair is not covered ($y^{rs}=0$). The definition of node coverage $z_r$ is provided in equation (8). Note that the variable $z_r$ is auxiliary variable for improved conciseness in the objective. The elimination of variable $z_r$ and replacement of variable $y^{rs}$ in the objective will not change the property of the formulation. Budget constraint is included in constraint (9). Constraints (10) are the boundary conditions of the initial state-of-refueling for all vehicles. Binary and nonnegativity constraints are included in inequalities (11)-(15). 

This MILP is essentially a maximum-coverage location problem \citep{daskin2011network}, which is NP-hard. Thus, we develop heuristic solutions based on genetic algorithm for solution efficiency. 
\section{A heuristic based on genetic algorithm}
\noindent 
There are various solution methods eligible for solving this MILP exactly or heuristically. The typical exact solution method is the branch-and-bound approach \citep{bertsimas1997introduction}. Though well suited for solving MILP, it often suffers the curse of dimensionality. Thus, many prior research efforts have been focused on the development of heuristic solution methods for solving facility location problems in a real-world context while balancing solution efficiency and quality \citep{melo2009facility}. One of the methods that have been widely used is the Lagrangian relaxation based heuristic solutions with applications in solving set-covering problems \citep{beasley1990lagrangian} and other location problems \citep{beasley1993lagrangean}. However, the success of implementation of the Lagrangian relaxation method is problem specific and depends on several factors, including the identification of constraint(s) to be relaxed, the goodness of bounds, and the solution efficiency of the relaxed problem. Other heuristics that have been proven to be effective for location problems include local search  \citep{korupolu2000analysis,ghosh2003neighborhood,arya2004local}, simulated annealing algorithm \citep{heragu1992experimental,murray1996applying,jayaraman2003simulated}, and genetic algorithm \citep{beasley1996genetic,houck1996comparison}.\\
\indent In this study, we adopt a heuristic based on genetic algorithm (GA) developed by \cite{beasley1996genetic}. Modifications are made to customize the procedures suitable for our particular problem. For completeness, the major procedures are reported and explained.\\
\indent \textbf{\emph{Representation and fitness function}}: A 0-1 $X_{|N|}$ array is used to represent location decisions over the transportation network, in which each bit of the array is notated as $X_i$ to indicate whether or not to place a refueling station at a particular node of the network. For example $X_i=1$ indicates a refueling station is built at node $i$. The array $X_{|N|}$ is called a feasible solution, and with the solution of location decision variables route decisions $Y^{rs,k}$ and detailed refueling decisions $B_i^{rs,k},l_i^{rs,k}$ are obtained by solving the original problem with the location decision variables fixed. The objective value is used as the the fitness value associated with this solution $X_{|N|}$.\\
\indent \textbf{\emph{Parent selection}}: In this study, we adopt the binary tournament parent selection and the fitness-based crossover method used in \citep{beasley1996genetic}. Initialized with a large population of solutions (e.g., 100), we randomly pick four individuals, to form two pools, each of which contains two individuals.\\
\indent \textbf{\emph{Crossover and mutation}}: The individual solution with better fitness in each pool is selected as a parent solution to breed a child solution based on fitness-based crossover operator \citep{beasley1996genetic}. Let $f_{P1}$ and $f_{P2}$ be the fitness values of parent $P1$ and $P2$, respectively and let $C$ be a child solution. 
\begin{enumerate}
\item[1)] if $P1 = P2$, then set $C:=P1$ or $P2$;
\item[2)] if $P1 \neq P2$, then set $C:=P1$ with probability $\dfrac{f_{P2}}{f_{P1}+f_{P2}}$, and $C:=P2$ with probability $\dfrac{f_{P1}}{f_{P1}+f_{P2}}$.
\end{enumerate}
\indent In addition to generating child solutions by crossover, child solution is also constructed by selecting the worst solution from the population and flipping the value for one of the bits of the solution according to the relationship of a random value and the mutation rate $p$. For example, if the random value generated within 0 and 1 is less than $p$, then the value of the bit is flipped.\\
\indent \textbf{\emph{Population replacement}}: If a child solution is identical to any of the solutions in the population, this child solution will be neglected; otherwise, it replaces the solution in the initial population with the worst fitness. \\
\indent Finally, we include complete and detailed procedures of the heuristic as follows:
\begin{enumerate}
\item Initialization: set values for parameters $m$, $n$, and $p$
\item While $i \leq n$ do 
\begin{enumerate}
\item Repeat following steps to generate $m$ child solutions
\begin{itemize}
\item Pick four solutions from the population to form two pools, each of which contains two solutions;
\item Select the solution with better fitness in each pool as one of the parents then do crossover to form a child solution;
\item Mutation is applied to the worst solution in the population with  mutation rate $p$;
\item Fix location decisions then solve the original problem by CPLEX to get the fitness value for the corresponding child solutions.
\end{itemize}
\item Replace the worst $m$ solutions in the population with child solutions generated in the iteration.
\end{enumerate}
\item Return the best solution in the population as the final solution to the original problem.
\end{enumerate}

\section{Numerical experiments}
\begin{figure}[ht]
\centering
\includegraphics[width=0.4\textwidth]{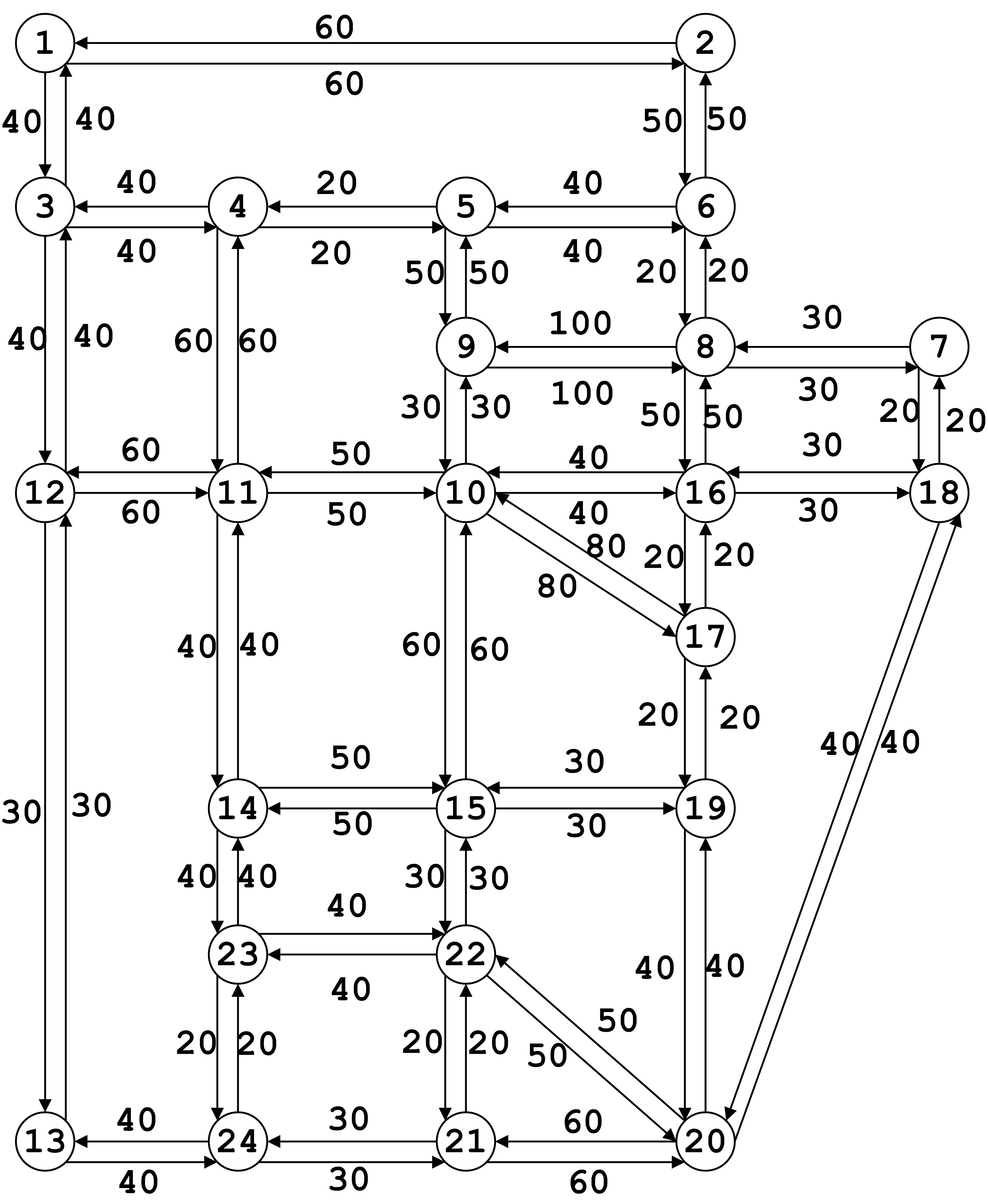}
\caption{The Sioux Falls network}
\label{SiouxFallsNetwork}
\end{figure}

\noindent We implemented the proposed P-MPRLM on the Sioux Falls network \citep{leblanc1975efficient}, which is a benchmark numerical example that has been widely used in transportation research community, to demonstrate the applicability of the model and gain insights about the effects of parameters on the model solutions. The Sioux Falls network shown in Figure \ref{SiouxFallsNetwork} is consisted of 24 nodes and 76 links with distance (assumed to be in miles) labeled on the on the links. In the numerical experiments, all nodes are both origins and destinations and candidate locations for AFSs. We assume that the probability of a node becoming a demand node is known a priori, which can be estimated based on statistical analysis (e.g., a result of logistic regression analysis \citep{hosmer2000introduction} on multivariate factors). In this study, the probability is randomly generated for the purpose of numerical experiments (see Table \ref{table:EVAdoptiongProbability}).

For each O-D pair, we consider multiple deviation paths between O-D pairs. The driving range of vehicles is set to be 100 miles. We solved the model both exactly and heuristically. For exact solutions, we first programmed the model by using AMPL and solved it to optimality by using CPLEX 12.6. For heuristic solutions, the model was programmed in MATLAB and solved by using our developed heuristic method based on genetic algorithm. Both exact and heuristic solutions were executed on a desktop with 8GB RAM and Intel Core i5-2500@3.30GHz processor under Windows 7 environment. In addition, all deviation paths are pre-generated exogenously using MATLAB and the solution times are not reported herein. The details on generating the deviation paths can be referred to \citep{li2014heuristic}. 

\begin{table}[ht]
\footnotesize
\centering
\caption{Probability of a node being an EV adopter}
\begin{tabular}{cccc}
\toprule
Node ID & Probability & Node ID & Probability\\
\hline
1 & 0.7689 & 13 & 0.7900\\
2 & 0.1673 & 14 & 0.3185\\
3 & 0.8620 & 15 & 0.5431\\
4 & 0.9899 & 16 & 0.0900\\
5 & 0.5144 & 17 & 0.1117\\
6 & 0.8843 & 18 & 0.1363\\
7 & 0.5880 & 19 & 0.6787\\
8 & 0.1548 & 20 & 0.4952\\
9 & 0.1999 & 21 & 0.1897\\
10 & 0.4070 & 22 & 0.4950\\
11 & 0.7487 & 23 & 0.1476\\
12 & 0.8256 & 24 & 0.0550\\
\bottomrule	
\end{tabular}
\label{table:EVAdoptiongProbability}
\end{table}

\subsection{Baseline case}
\noindent We implement the P-MPRLM with up to three deviation paths (K=3) considered for each O-D pair. We considered 12 different budget level (i.e., in a range of one to 12 stations). Note that although all 24 nodes are candidate AFS locations, 12 stations are sufficient to fully cover all the nodes (i.e., $z_r=1, \forall r \in R$). Here, with $C_i, \forall i \in \tilde{N}$ being unity, the budget $m$ is essentially the total number of AFSs to be deployed on the network.

We first evaluate the \emph{quality} and \emph{efficiency} of the heuristic solution by comparing with the counterparts of CPLEX solutions. Table \ref{table: Comparison of objective values between Heuristic and CPLEX} reports the objective values (i.e., the total expected node coverage) of the heuristic and CPLEX solutions and their difference under 12 different budget levels. For each budget level, we run the heuristic for 50 times and report the average and CPLEX solutions are the maximum total expected node coverage of all nodes. From the table, the differences in objective values between CPLEX and our heuristic are in a small range between 0.0\% and 1.9\%, which well justifies the quality of heuristic solutions. We then use the results of budget of three stations in Table \ref{table: Details of coverage when budget is available for built at most 3 stations} to explain what the sum of expected node coverage is composed of. First, under the budget, three stations are deployed at nodes \# 3, \#6 and \#16 respectively. The details on the coverage of OD pairs for each node are presented in the table and the total expected node coverage is on the last row. 

\begin{figure}[ht]
\centering
\begin{subfigure}{0.9\textwidth}
  \centering
  \includegraphics[width=1\linewidth]{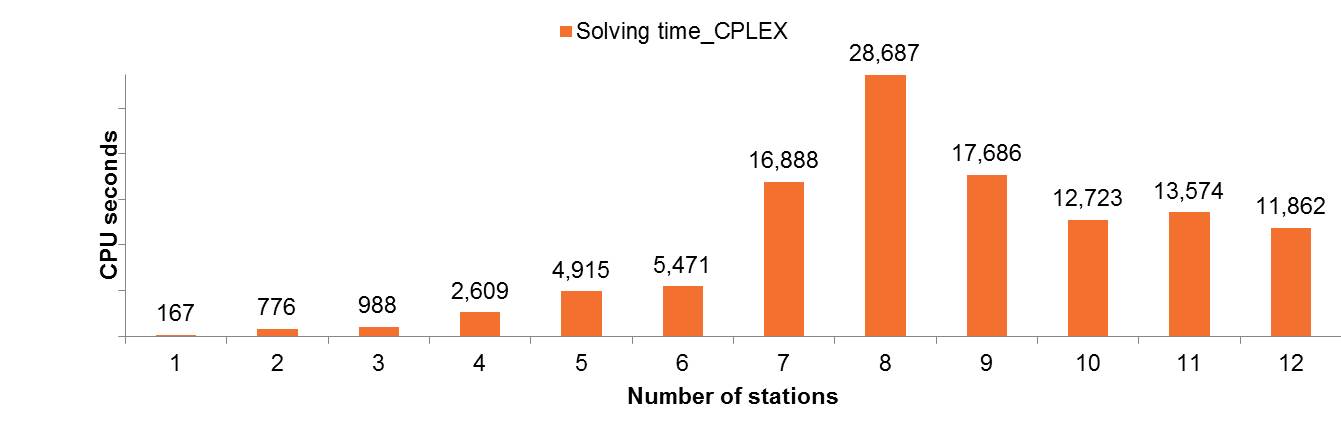}
  \caption{Solution time by CPLEX}
  \label{SolvingTime:Cplex}
\end{subfigure}
\begin{subfigure}{0.9\textwidth}
  \centering
  \includegraphics[width=1\linewidth]{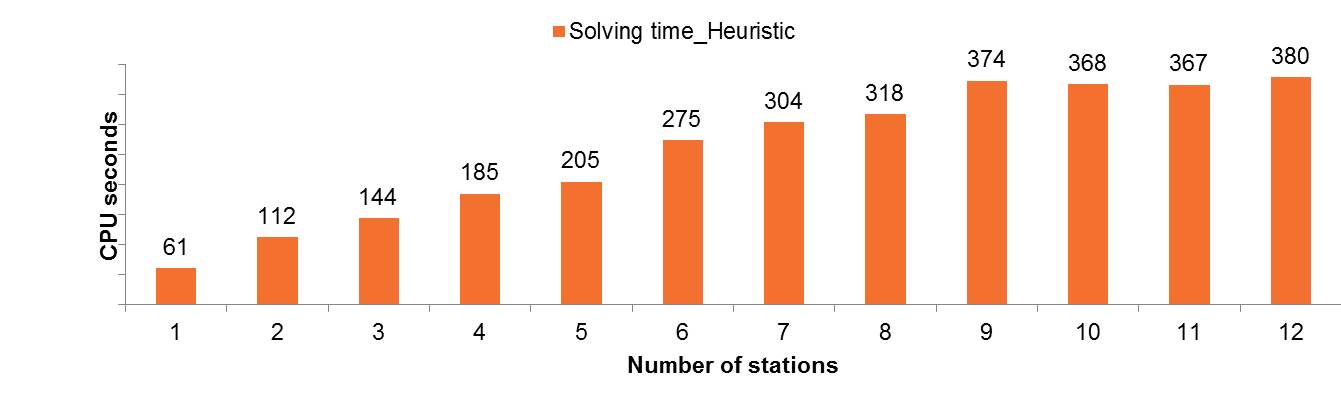}
  \caption{Solution time by the heuristic}
  \label{SolvingTime:Heuristic}
\end{subfigure}
\caption{Computational times by Heuristic and CPLEX}
\label{SolvingTime}
\end{figure}

Now, let us compare the solution efficiency in terms of CPU seconds. The average solution times of 50 runs of the heuristic for each of the 12 budget levels are reported in Figure \ref{SolvingTime:Heuristic}, which is compared with the solution times of counterparts solved by CPLEX in Figure \ref{SolvingTime:Cplex}. From the results, heuristic significantly reduces the solution times compared to CPLEX in all cases. The largest difference is when eight stations to be deployed and the solution time by CPLEX is about 90 times longer than the heuristics. Also, we note that the problem becomes more substantially complicated with the increase of budget (i.e., more candidate stations to be considered). From the results of the baseline numerical experiments, we are confident that the heuristic is solution efficient without compromising high solution quality.

\begin{table}[ht]
	\footnotesize
	\centering
	\caption{Comparison of objective values between Heuristic and CPLEX}
	\begin{tabular}{cccc}
		\toprule
		Budget & Objective value & Objective value& Difference\\
		(Number of stations)  & Heuristic & CPLEX & (\%)\\
		\hline
		1 &  2.45 & 2.45 & 0.00\% \\
		2 & 3.78 & 3.79 & 0.20\%\\
		3 & 5.04 & 5.11 & 1.40\%\\
		4 & 6.29 & 6.36 & 1.10\%\\
		5 & 7.41 & 7.54 & 1.70\%\\
		6 & 8.49 & 8.58 & 1.10\%\\
		7 & 9.20 & 9.29 & 1.00\%\\
		8 & 9.69 & 9.88 & 1.90\%\\
		9 & 10.22 & 10.33 & 1.10\%\\
		10 & 10.45 & 10.52 & 0.70\%\\
		11 & 10.62 & 10.66 & 0.40\%\\
		12 & 10.67 & 10.69 & 0.20\%\\
		\bottomrule	
	\end{tabular}
	\label{table: Comparison of objective values between Heuristic and CPLEX}
\end{table}

\begin{table}[ht]
	\footnotesize
	\centering
	\caption{Details of coverage when budget is available for built at most 3 stations}
	\begin{tabular}{cccc}
		\toprule
		Node ID & Probability ($p_r$) & Node coverage ($z_r$) & $p_rz_r$\\
		\hline
		1 &  0.7689 & 0.42 & 0.32 \\
		2 & 0.1673 & 0.54 & 0.09\\
		3 & 0.8620 & 0.46 & 0.40\\
		4 & 0.9899 & 0.54 & 0.54\\
		5 & 0.5144 & 0.58 & 0.30\\
		6 & 0.8843 & 0.50 & 0.44\\
		7 & 0.5880 & 0.50 & 0.29\\
		8 & 0.1548 & 0.50 & 0.08\\
		9 & 0.1999 & 0.42 & 0.08\\
		10 & 0.4070 & 0.46 & 0.19\\
		11 & 0.7487 & 0.54 & 0.41\\
		12 & 0.8256 & 0.38 & 0.31\\
		13 & 0.7900 & 0.33 & 0.26\\
		14 & 0.3185 & 0.13 & 0.04\\
		15 & 0.5341 & 0.54 & 0.29\\
		16 & 0.0900 & 0.50 & 0.05\\
		17 & 0.1117 & 0.50 & 0.06\\
		18 & 0.1363 & 0.54 & 0.07\\
		19 & 0.6787 & 0.54 & 0.37\\
		20 & 0.4952 & 0.46 & 0.23\\
		21 & 0.1897 & 0.17 & 0.03\\
		22 & 0.4950 & 0.50 & 0.25\\
		23 & 0.1476 & 0.17 & 0.02\\
		24 & 0.0500 & 0.17 & 0.01\\
		\hline
		Expected node coverage & & & 5.11\\
		\bottomrule	
	\end{tabular}
	\label{table: Details of coverage when budget is available for built at most 3 stations}
\end{table}

\begin{figure}[ht]
\centering
\includegraphics[width=0.45\textwidth]{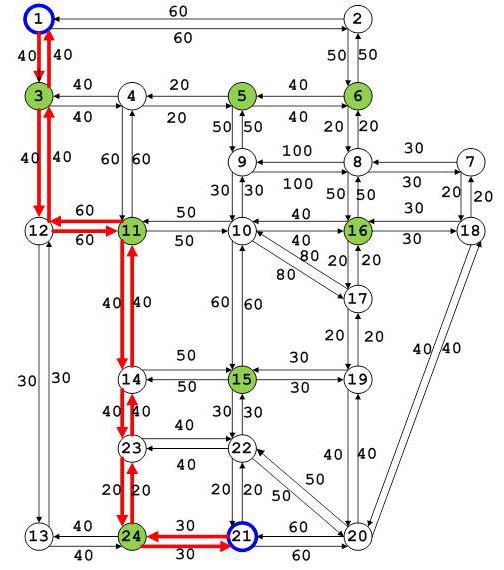}
\caption{Location and routing decisions when VR=100 and K=3 with seven stations}
\label{LocationAndRoute}
\end{figure}

\indent We take the scenario of locating seven stations on the network (see Figure \ref{LocationAndRoute}) as an example to demonstrate the integrated decisions on station locations, route choices, and refueling schedules. In the figure, selected stations are highlighted. Note that the location pattern of stations may  not be unique, but the resulting expected node coverage is maximum. Let us use the O-D pair (1, 21) as an example to show how routing and refueling are compatible with station location decisions. The round-trip between the O-D pair is completed via the second shortest path $1\leftrightarrow3\leftrightarrow12\leftrightarrow11\leftrightarrow14\leftrightarrow23\leftrightarrow24\leftrightarrow21$ , as highlighted in the figure. 

\begin{figure}[ht]
\centering
\includegraphics[width=0.8\textwidth]{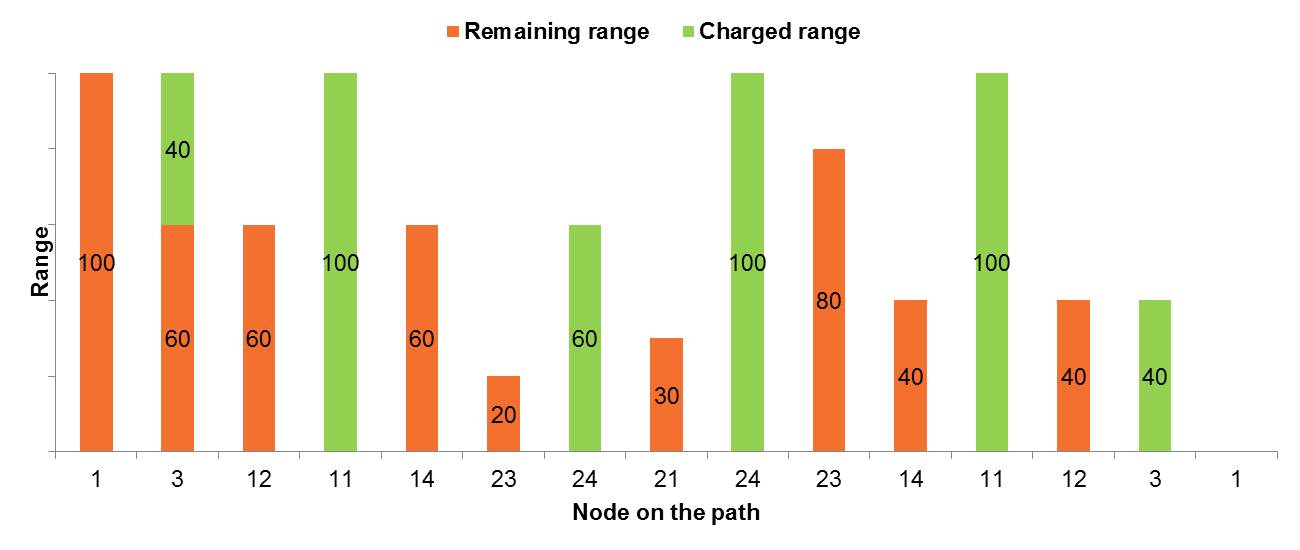}
\caption{Detailed refueling schedule for O-D pair (1,21)}
\label{ChargingSchedule}
\end{figure}

The detailed refueling schedule along the path is shown in Figure \ref{ChargingSchedule}, in which the orange bars are the remaining on-board fuel before refueling at a particular node while the green bars are the range extended by refueling at stations. A combined bar denotes the state-of-fueling in terms of driving range when departing from the node. From the figure, multiple refueling stops are needed to complete the round trip between the node pair (1, 21). 

\subsection{Sensitivity Analysis}
We conducted a series of sensitivity analyses to understand the impacts of vehicle range, initial state of fuel and probabilities on the decisions on station locations and resulting expected node coverage. 

\begin{figure}[ht]
\centering
\includegraphics[width=0.7\textwidth]{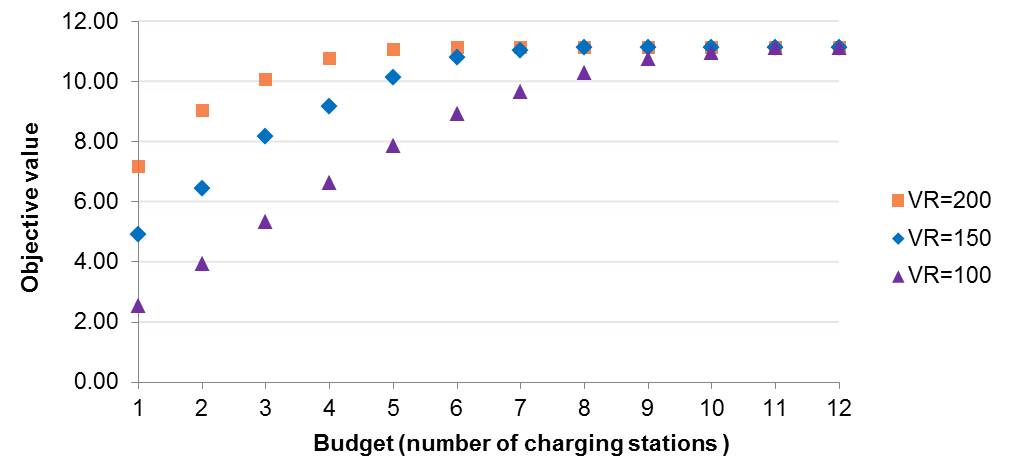}
\caption{Impact of budget and vehicle range}
\label{BudgetAndVR}
\end{figure}

\textbf{Vehicle Range (VR).} We consider three vehicle ranges (i.e., VR=100, 150 and 200) and run the model under all 12 budget levels and three deviation paths (K=3). The resulting total expected node coverage (i.e., objective values) are reported in Figure \ref{BudgetAndVR}. From the results, an extended vehicle range will help improve node coverage at every given budget level and the effect is more substantial when the budget level is relatively lower. For example, when only one station is placed, the total expected node coverage when vehicle range is 200 miles is 7.19, which almost three times as of the coverage of 2.55 when the vehicle range is 100 miles. The effect becomes marginal at a raised budget (e.g., nine stations or more) and the difference in terms of the expected node coverage is indiscernible. This is because with more stations on the network, drivers have increased accessibility to refueling and thus the restraining effects on vehicle range become weaker. 

For a particular vehicle range, there may exist a critical number of stations, by which the maximum node coverage has already reached. In other words, additional available stations will increase the redundancy of refueling opportunities, but cannot further improve the demand node coverage. For example, when vehicle range is 100 miles, 11 stations is identified as the critical number of stations. Similarly, the such numbers for vehicle ranges of 150 and 200 miles are seven and five stations respectively. The critical number of stations decreases with the increase of vehicle range. This is because as vehicle range is extended, fewer refueling stations are needed for achieving the same node coverage. 
From the results, we may immediately identify the maximum coverage for a given vehicle range and budget level and the minimum budget for the maximum coverage given a vehicle range. 

\begin{figure}[ht]
	\centering
	\includegraphics[width=0.7\textwidth]{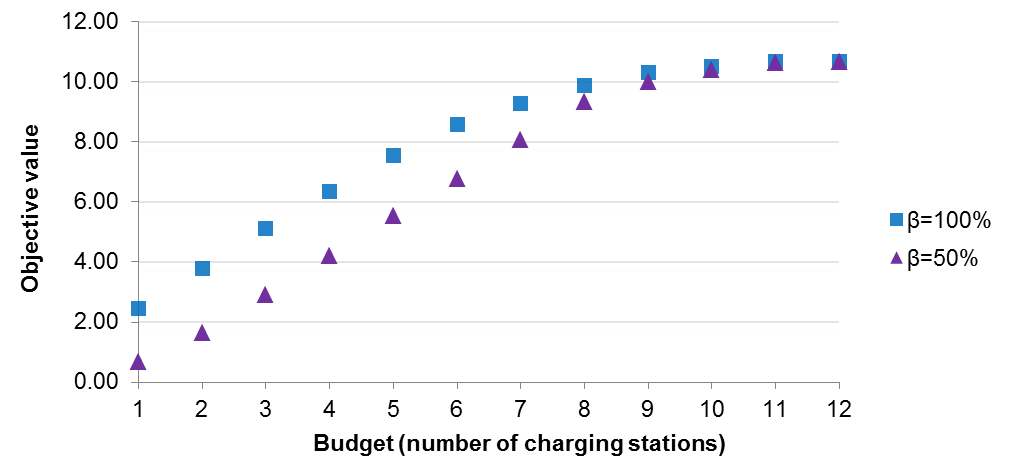}
	\caption{Impact of initial battery status}
	\label{BatteryStatus}
\end{figure}

\textbf{Initial state-of-fueling (SOF).} This set of analyses on the initial SOF aims to gain insights about how the total node coverage will be impacted as a result of the initial SOF. A great deal of prior research efforts have emphasized on the assumption of initial SOF or fuel status, which is believed to be substantially impactful on the optimal solutions \citep{kuby2005flow}. In our baseline case, the initial SOF, denoted by $\beta$ in the model, is assumed to be 100\%, meaning fully refueled at origins. In this sensitivity analysis, we tested a case when the initial SOF is assumed to be half full (i.e., $\beta=50\%$VR) for all the nodes and when $VR=100$ miles, which is a common assumption held in prior research studies \citep{kuby2005flow} when round-trips were considered. The resulting node coverage is reported in Figure \ref{BatteryStatus}, compared with the baseline results (i.e., $\beta=100\%$VR). It is obvious that for a given budget, lower initial SOF results in lower expected node coverage. As more stations become available, the difference in coverage (i.e., objective value) decreases. When there are sufficient number of stations (e.g., 10 stations), the difference in node coverage is negligible, which is less than 1\%. The reason can be similarly explained. Increasingly available stations make refueling more accessible and thus both the initial fueling status and vehicle range become less crucial.

\begin{figure}[ht]
	\centering
	\includegraphics[width=0.7\textwidth]{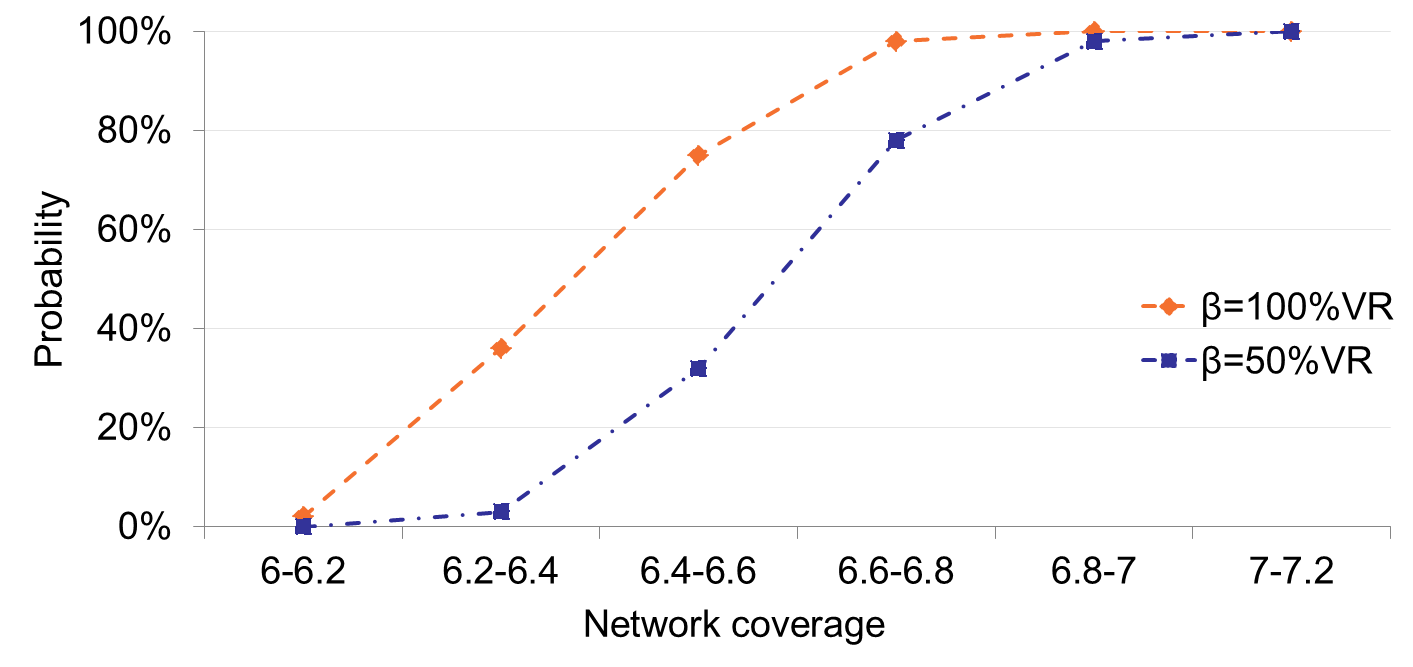}
	\caption{Impact of full and half initial SOFs}
	\label{CDFfullandhalfSOFs}
\end{figure}

In reality, however, AFV users may have different initial SOFs. To investigate the performance of the deployments that are result of constant initial SOFs ($\beta=100\%$VR and $\beta=50\%$VR), we conducted another numerical experiment by randomly generating 100 initial SOFs between 0 and 100\% for each budget level. We then evaluated the solutions of $\beta=100\%$VR and $\beta=50\%$VR in terms of the demand node coverage and report the average of the sum of the expected node coverage (i.e., averaged objective value over 100 scenarios) under each budget level in Table \ref{AverageCoverageOfRandomSOFs}. From the results, the assumption of $\beta=50\%$VR overall outperformed the assumption of $\beta=100\%$VR, indicated by the overall higher average sum of expected node coverage, which is understandable as $\beta=50\%$VR is a more conservative assumption of initial SOF than the $\beta=100\%$VR. In fact, the differences between the assumptions are trivial for all different budget levels, which implies that the assumptions of initial SOFs in the model may not be as sensitive to the effects of the resulting location planning strategies on the Sioux Falls network. However, note that this is the average of resulting node coverage and it can highly be dependent on the topological structure of the network. We further investigate the full spectrum of the distributions of the sum of node coverage for the 100 scenarios and report them in a cumulative probability distribution against the sum of node coverage for a given budget of seven stations, in Figure\ref{CDFfullandhalfSOFs}. The figure presents a clear shift of the curve to the right when the initial SOF is reduced from 100\%VR TO 50\%VR. For example, for the same 40\% of chance, the assumption of $\beta=50\%$VR will result in about a coverage of 6.6, but the assumption of $\beta=100\%$VR will lead to a lower coverage of about 6.4. This observation is consistent with the results of average node coverage that a more conservative assumption of $\beta=50\%$VR will help improve the node coverage on the network slightly better than the assumption of full initial SOF. Eventually, both assumptions will result in the same highest coverage at 7.2.

\begin{table}[htbp]
	\centering
	\caption{The average sum of expected node coverage under 100 random initial SOFs}
	\begin{tabular}{ccc}
		\toprule
		Budget (Number of stations) & $\beta=100\%$VR & $\beta=50\%$VR\\
		\midrule
		1    & 0.74 & 0.74 \\
		2    & 1.49 & 1.45\\
		3    & 2.60 & 2.60\\
		4    & 3.62 & 3.52\\
		5    & 4.34 & 4.75\\
		6     & 5.55 & 5.62\\
		7     & 6.48 & 6.67 \\
		8    & 7.62 & 7.62\\
		9     & 8.01 & 8.14\\
		10    & 8.48 & 8.68\\
		11    & 8.90 & 8.94\\
		12     & 8.98 & 9.13\\
		\bottomrule
	\end{tabular}%
	\label{AverageCoverageOfRandomSOFs}%
\end{table}%

\begin{figure}[ht]
\centering
\includegraphics[width=0.7\textwidth]{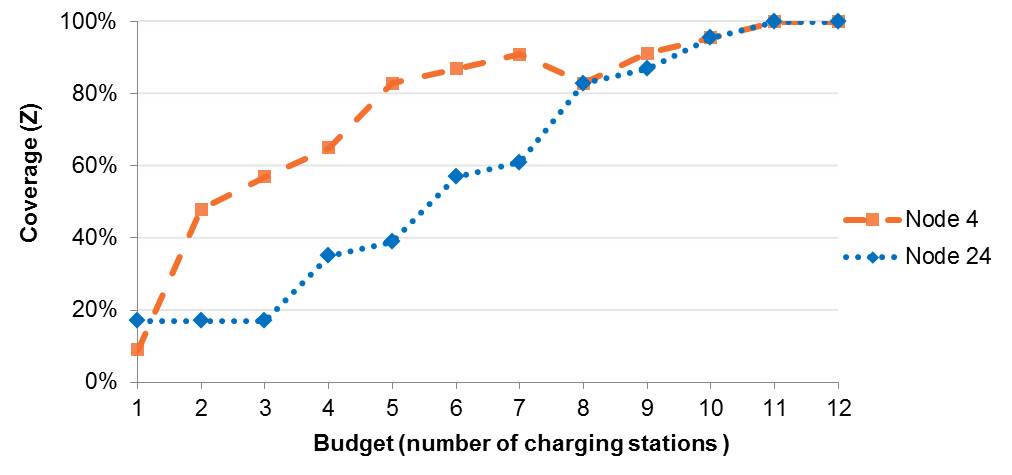}
\caption{Coverage ($z$) for node $\#4$ and $\#24$ when VR=100 and K=3}
\label{IndividualCityCoverage}
\end{figure}

\textbf{The effects of probability.} 
The decisions of stations on the network and resulting node coverage are the combined effects of probability of demand nodes and their geographic locations on the network. We are particularly interested in understanding the potential trade-off between the probability and geographic location of a node and their impacts on the expected node coverage. To our interest, we pick two nodes on the network to analyze, which are associated with the highest ($\#4$) and lowest ($\#24$) probabilities (i.e., 0.9899 and 0.0550), respectively. Their expected coverage under different budget levels are plotted in Figure \ref{IndividualCityCoverage}. Node \#4, though with highest probability, is not always guaranteed to receive higher coverage than the node \#24. In fact, when budget is extremely low (e.g., only one station placed), node \#24 has higher coverage than node \#4. To understand this phenomenon, we need to refer to their geographic distributions on the network in Figure \ref{SiouxFallsNetwork}. 
The location of node \#24 dominates the coverage as trips between some of O-D pairs (e.g., (24, 13), (24, 21), and (24, 23)) can be self completed without the need of refueling while the node \#4 though having highest probability has less privileged geographic location compared to node \#24 in the network, given the vehicle range of 100 miles. However, with the increase of budget levels, as node \#4 has significantly higher probability than node \#24, the expected coverage of node \#4 rises rapidly. For example, the node coverage increases from 9\% to 57\% when the number of stations increases from 1 to 3, while the expected coverage of node \#24 is much lower. However, when the network is present with sufficiently large number of refueling stations, the expected coverage for both nodes is almost equivalent. 


It is also noted that coverage of node \#4 drops as number of stations increases from seven to eight, while at the same time coverage of node \#24 increases. This is a perfect demonstration of system optimization. Some individual node (e.g., \#4) coverage is sacrificed in order to improve the overall expected node coverage for the entire network.

We are also interested in understanding the effects of incorporating the probability of demand nodes in the decision of siting stations on the network through comparing the counterparts of not having probability involved, i.e., set $p_r=1$. We conducted numerical experiments on a particular budget level of seven stations. We only report the coverage of nodes ($z_r$) with probabilities either greater than 0.8 (four nodes) or lower than 0.2 (nine nodes) when VR=100 and K=3 in Table \ref{ComparisonOfCoverageBetweentWithAndWithoutProbability} . The coverage of nodes without incorporating probability is a result of stripping out the effects of probability and is only driven by the geographic distributions of nodes on the network. By comparing the coverage of without probability to the counterparts of having probability involved, all nine low-probability nodes receive as least as much coverage with an increase in a range between 0\% and 8\% while three out of four high-probability nodes obtain as most as much coverage with a decrease in a range between 0\% and 5\%. From the results, probability needs to be incorporated to counter the effects of the geographic locations in order to achieve the maximum node coverage on a system level.





\begin{table}[htbp]
  \centering
  \caption{Comparison of coverage of nodes with high( $\geq0.80$) and low ($\leq0.20$) probabilities}
    \begin{tabular}{cccc}
    \toprule
    Node ID & Probability of demand nodes & Coverage\_probability & Coverage\_no probability \\
    \midrule
    24    & 0.0550 & 61\%  & 78\% \\
    16    & 0.0900 & 87\%  & 87\% \\
    17    & 0.1117 & 83\%  & 87\% \\
    18    & 0.1363 & 83\%  & 91\% \\
    23    & 0.1476 & 74\%  & 78\% \\
    8     & 0.1548 & 87\%  & 87\% \\
    2     & 0.1673 & 78\%  & 83\% \\
    21    & 0.1897 & 74\%  & 87\% \\
    9     & 0.1999 & 87\%  & 87\% \\
    \midrule
    12    & 0.8256 & 96\%  & 91\% \\
    3     & 0.8620 & 87\%  & 87\% \\
    6     & 0.8843 & 83\%  & 87\% \\
    4     & 0.9899 & 91\%  & 87\% \\
    \bottomrule
    \end{tabular}%
  \label{ComparisonOfCoverageBetweentWithAndWithoutProbability}%
\end{table}%

\clearpage
\section{Conclusion and Future Work}
\noindent This paper presents a probability-based multi-path refueling location model that incorporates AFV adoption probabilities of demand cities, vehicle range, users deviation choice and round-trip between O-D pairs to maximized the total weighted coverage of all demand cities with a given budget. A genetic algorithm based heuristic is adopted to help efficiently solve the proposed problem. Numerical experiments are implemented with the Sioux Falls network to demonstrate the model and heuristic, and sensitivities analysis are conducted to examine the impacts of vehicle range, initial state of fuel, and AFV adoption probability. The results indicate that by integrating probability information in the model, cities with high and low potential of adopting EVs would experience improvement and decreases in coverage, which would provide insights to better utilize the limited resources. \\
\indent Future research can be conducted in a few directions. First, in this paper the probability of each city is fixed, while in reality it may subject to uncertainty due to socio-economic developments. Thus, more accurate prediction method or discrete scenarios can be considered to further improve this study. Second, refueling time is not considered in this paper and in the future studies can be conducted to investigate how it could affect users' route choice and station location decisions. In further, once AFVs (e.g., EVs) are massively adopted, AFV traffic flows can be integrated and congestion effects at refueling stations may become more significant, and in that case smart capacity design and pricing scheme may be considered.


\end{document}